\documentclass[a4paper,12pt,reqno]{amsart}

\usepackage[utf8]{inputenc}
\usepackage[T1]{fontenc}
\usepackage{lmodern}
\usepackage[english]{babel}
\usepackage{microtype}

\usepackage{amsmath,amssymb,amsfonts,amsthm,comment}
\usepackage{mathtools,accents}
\usepackage{mathrsfs}
\usepackage{aliascnt}
\usepackage{braket}
\usepackage{bm}
\usepackage{esint}
\usepackage{amsfonts}
\usepackage{csquotes}

\usepackage[a4paper,twoside,top=0.9in, bottom=0.9in, left=0.85in, right=0.85in]{geometry} 
\usepackage{hyperref}
\usepackage[english,capitalize]{cleveref}

\usepackage{enumerate}
\usepackage{xcolor}
\usepackage{comment}

\usepackage{imakeidx} 
\usepackage{xparse} 
\usepackage{mathtools}
\usepackage{float}
\usepackage{comment}
\usepackage{nicefrac}
\usepackage{graphicx}
\usepackage{enumerate}
\usepackage{enumitem}

\definecolor{newblue}{rgb}{0.2, 0.3, 0.85}
\hypersetup{colorlinks=true, linkcolor=newblue, citecolor=newblue, urlcolor  = newblue}

\usepackage[maxbibnames=99,backend=biber, sorting=nyt, doi=false, url=false, isbn=false]{biblatex}
\makeatletter
\g@addto@macro\@floatboxreset\centering
\makeatother

\makeatletter
\def\newaliasedtheorem#1[#2]#3{
  \newaliascnt{#1@alt}{#2}
  \newtheorem{#1}[#1@alt]{#3}
  \expandafter\newcommand\csname #1@altname\endcsname{#3}
}
\makeatother

\numberwithin{equation}{section}

\newtheoremstyle{slanted}{\topsep}{\topsep}{\slshape}{}{\bfseries}{.}{.5em}{}

\theoremstyle{plain}
\newtheorem{theorem}{Theorem}[section]
\newaliasedtheorem{proposition}[theorem]{Proposition}
\newaliasedtheorem{question}[theorem]{Question}
\newaliasedtheorem{problem}[theorem]{Problem}
\newaliasedtheorem{lemma}[theorem]{Lemma}
\newaliasedtheorem{corollary}[theorem]{Corollary}
\newaliasedtheorem{counterexample}[theorem]{Counterexample}
\newtheorem*{acknowledgement}{Acknowledgement}

\theoremstyle{definition}
\newaliasedtheorem{definition}[theorem]{Definition}
\newaliasedtheorem{example}[theorem]{Example}
\newaliasedtheorem{conjecture}[theorem]{Conjecture}

\theoremstyle{remark}
\newaliasedtheorem{remark}[theorem]{Remark}

\newcommand{\N}{\mathbb{N}}

\newcommand{\R}{\mathbb{R}}

\newcommand{\nchi}{{\raise.3ex\hbox{\(\chi\)}}}

\let\phi\varphi

\newcommand{\abs}[1]{\left\lvert#1\right\rvert}
\newcommand{\norm}[1]{\left\lVert#1\right\rVert}

\DeclareMathOperator{\sign}{sign}

\newcommand{\dist}{\mathsf{d}}

\newfont{\tmpf}{cmsy10 scaled 2500}

\newcommand{\de}{\ensuremath{\,\mathrm d}} 

\author[Antonelli]{Gioacchino Antonelli}
\address{Courant Institute Of Mathematical Sciences (NYU), 251 Mercer Street, 10012, New York, USA}
\email{ga2434@nyu.edu}

\author[Calzi]{Mattia Calzi}
\address{Dipartimento di Matematica Federigo Enriques\\
Università degli Studi di Milano\\
Via Saldini, 50
20133 Milano (MI), Italy}
\email{mattia.calzi@unimi.it}

\author[Gordina]{Maria Gordina
}
\address{ Department of Mathematics\\
University of Connecticut\\
Storrs, CT 06269,  U.S.A.}
\email{maria.gordina@uconn.edu}


\begin{document}

\title[Logarithmic Sobolev inequalities]{Sharp defective log-Sobolev inequalities on $H$-type groups}

\begin{abstract}
In this paper we prove  a sharp defective log-Sobolev inequality on $H$-type groups. Then we use such an inequality to show exponential integrability of Lipschitz functions with respect to the heat kernel measure. A defective log-Sobolev-type inequality  for the Gaussian-like measure with respect to the sub-Riemannian distance is also proved on arbitrary $H$-type groups.
\end{abstract}

\maketitle

\section{Introduction}

This paper deals with logarithmic Sobolev inequalities on Lie groups. This line of research was pioneered by L.~Gross in \cite{GrossLSIGroups} and has recently seen many developments in terms of settings and techniques. In particular, we are motivated by the following (open) question.

Let $\mathbb H^n$ be the (isotropic) $n$th Heisenberg group. Let $\de\mu$ be a Haar measure on $\mathbb H^n$, let $\Delta_{\mathbb H^n}$ be the sub-Laplacian associated to the restriction of $\langle\cdot,\cdot\rangle$ to the horizontal distribution of $\mathbb H^n$, and denote by $\de\mu_t$ the heat kernel measure corresponding to the semigroup $(e^{t\Delta_{\mathbb H^n}})_{t\geqslant 0}$. For details on such measures we refer to  \cite[p.~951]{DriverGrossSaloff-Coste2009a}. 
It is an \emph{open problem} to determine the best constant $\omega>0$ (and find optimizers, if they exist) for which for every $t>0$, and every $f\in C^{\infty}_c(\mathbb H^n)$ with $\int_{\mathbb H^n} f^2\, \de\mu_t =1$, one has
\begin{equation}\label{eqn:LogSobolevClassica}
\int_{\mathbb H^n} f^2\log(f^2)\,\de\mu_t\leqslant \omega t\int_{\mathbb H^n} |\nabla_{\mathbb H^n} f|^2\,\de\mu_t,
\end{equation}
where $|\nabla_{\mathbb H^n}u|$ is the norm of the horizontal gradient $\nabla_{\mathbb H^n}u$ of $u$. The existence of a finite positive constant $\omega$ for which \eqref{eqn:LogSobolevClassica} holds is due to \cite{DriverMelcher2005, LiHQ2006, LiHQ2007}, and the discussion in \cref{sec:logsobolevCD}. In the hope to shed a new light on \eqref{eqn:LogSobolevClassica}, this note deals with a sharp defective version of \eqref{eqn:LogSobolevClassica}. In particular, we obtain a sharp defective log-Sobolev inequality with computable constants on any $H$-type group equipped with a sub-Riemannian structure, and derive some consequences.

\subsection{Statements of the results}
Before stating the main theorem, we need to introduce some basics on $H$-type groups. We refer the reader to \cref{HtypeGroups} for a more detailed discussion.

Let $\mathbb G$ be an $H$-type group, and let $\langle\cdot,\cdot\rangle$ be an inner product on the Lie algebra $\mathfrak{g}$ for which \cref{item1Htype}, and \cref{item2Htype} hold in \cref{HtypeGroups}. Recall that $\langle\cdot,\cdot\rangle$ on $\mathfrak{g}$ gives rise to a left-invariant Riemannian metric on $\mathbb G$. Denote by $2n$ the dimension of the first layer of $\mathbb G$, and denote by $m$ the dimension of the second layer of $\mathbb G$. Let $\de\mu$ be a Haar measure, and observe that $\de\mu$ is bi-invariant. Let $\Delta_{\mathbb G},\mathrm{d}\mu_t$ be the sub-Laplacian and the corresponding heat kernel measure defined as in the case $\mathbb G=\mathbb H^n$ above. Note that by \cite[Section~3]{DriverGrossSaloff-Coste2009a} the operator $\Delta_{\mathbb G}$ is essentially self-adjoint on $C_{c}^{\infty}\left( \mathbb G\right)$. 
Let $|\nabla_{\mathbb G}u|$ denote the norm of the horizontal gradient $\nabla_{\mathbb G}u$ of $u$.

We denote by $W^{1,2}(\mathbb G,\de\mu_t)$  the completion
of $C^{\infty}_c(\mathbb G)$ with respect to the norm
\[
\|u\|_{W^{1,2}(\mathbb G,\de\mu_t)}^2:=\int_{\mathbb G}u^2\de\mu_t+\int_{\mathbb G}|\nabla_{\mathbb G} u|^2\de\mu_t.
\]
Note that $W^{1,2}(\mathbb G,\de\mu_t)$ is the Sobolev space defined as the domain of the Dirichlet form induced by the horizontal gradient, as done in \cite[p.~951]{DriverGrossSaloff-Coste2009a}, and \cite[p.~1902]{CarfagniniGordina2024}.

\begin{definition}
    Let $\vartheta,\eta \geqslant 0$. We say that the $H$-type group $\mathbb G$
    \emph{satisfies the defective log-Sobolev inequality with constants $\vartheta,\eta$}, shortly $\mathrm{DLS}(\vartheta,\eta)$, if the following holds. For every $t>0$, and for every $f\in W^{1,2}(\mathbb G,\de\mu_t)$, we have
\begin{equation}\label{eqn:FullDefectiveLogSob}
\int_{\mathbb G} f^2\log(f^2)\de\mu_t-\left(\int_{\mathbb G} f^2\de\mu_t\right)\log\left(\int_{\mathbb G}f^2\de\mu_t\right)\leqslant \eta \int_{\mathbb G} f^2\de\mu_t+\vartheta t\int_{\mathbb G} |\nabla_{\mathbb G} f|^2  \de\mu_t.
\end{equation}
\end{definition}
\begin{remark}
By a simple normalization argument, $\mathbb G$ satisfies $\mathrm{DLS}(\vartheta,\eta)$ if and only if the following holds. For every $t>0$, and every $f\in W^{1,2}(\mathbb G,\de\mu_t)$ with $\int_{\mathbb G} f^2 \de\mu_t=1$, we have
\begin{equation}\label{eqn:LogSobolevTypeUpToAdditive}
\int_{\mathbb G} f^2\log (f^2) \de\mu_t\leqslant \eta+\vartheta t \int_{\mathbb G}|\nabla_{\mathbb G} f|^2\de\mu_t.
\end{equation}
\end{remark}

We are now ready to state the main result of the paper.
\begin{theorem}\label{thm:NonAbelianHtype}
Let $G$ be an $H$-type group. Define
\[
\mathrm{DLS}(\mathbb G):=\{\vartheta >0: \text{there is $\eta\geqslant 0$ such that $\mathbb G$ satisfies $\mathrm{DLS}(\vartheta,\eta)$}\}.
\]
Hence 
\[
\inf\mathrm{DLS}(\mathbb G) = 4,
\]
and 
\begin{equation}\label{eqn:4notin}
4\notin \mathrm{DLS}(\mathbb G).
\end{equation}
\end{theorem}
\begin{remark}
    For every $\vartheta>4$ a constant $\eta$ for which $\mathbb G$ satisfies $\mathrm{DLS}(\vartheta,\eta)$ is computed in \eqref{eqn:ExplicitEta}.
    We stress that such an $\eta$ only depends on $\vartheta,m,n$, where $2n$ is the topological dimension of the first layer of $\mathbb G$, and $m$ is the topological dimension of the second layer of $\mathbb G$.
\end{remark}
The proof of \cref{thm:NonAbelianHtype} is in \cref{proof2}. Note that  Theorem~\ref{thm:NonAbelianHtype} says that, up to an additive constant,  we can take $\vartheta$ as close as we want to $4$ in \eqref{eqn:LogSobolevTypeUpToAdditive}, but there is no choice of $\eta\geqslant 0$ for which we can choose $\vartheta = 4$. This is in sharp contrast with what happens on $\mathbb R^n$. Indeed, in $\mathbb R^n$ the sharp log-Sobolev inequality with respect to the Gaussian measure states that for every $f\in W^{1,2}(\mathbb R^n)$, and every $t>0$, we have
\[
\int_{\mathbb R^n} f^2\log(f^2)\de\mu_t-\left(\int f^2\de\mu_t\right)\log\left(\int_{\mathbb R^n}f^2\de\mu_t\right)\leqslant 
4 t\int_{\mathbb R^n} |\nabla f|^2  \de\mu_t.
\]
The contrast put forward by Theorem~\ref{thm:NonAbelianHtype} can be seen as an effect of the fact that the sub-Riemannian manifold $\mathbb G$ does not have the Ricci curvature bounded from below in the weak sense. Compare with the discussion in \cref{sec:logsobolevCD}.

The proof of \cref{thm:NonAbelianHtype} combines two ingredients. The first is a refined heat kernel estimate on $H$-type groups in  \cref{lem:UnboundedOmega4}, for which we crucially exploit the delicate estimates for all the derivatives of the heat kernel obtained in \cite{BrunoCalzi2018}. The second is a consequence of the $(2,2^{\ast})$-Sobolev inequality on $H$-type groups in  \cref{cor:BecknerPearson}.

To the authors' knowledge, the result in \cref{thm:NonAbelianHtype} is new even when $\mathbb G$ is the standard isotropic $n$th Heisenberg group $\mathbb H^n$.
Moreover, the result in \cref{thm:NonAbelianHtype} appears to be one of the few instances where the constants in a log-Sobolev-type inequality on homogeneous groups are explicit. Indeed, as of today, not only the best constant in \eqref{eqn:LogSobolevClassica} is unknown, but also no explicit upper bound on such constant is currently available. On the contrary a lower bound on the best constant for which \eqref{eqn:LogSobolevClassica} holds is known in $\mathbb H^n$ (\cite{DriverMelcher2005, LiHQ2006}), and in $H$-type groups \cite{EldredgePhDThesis}. See \cref{thm:NonAbelianHtypePre} for additional details.
\medskip

It is well-known, leveraging Herbst's argument, that a log-Sobolev inequality implies concentration inequalities of Gaussian type. This remains true even when a defective log-Sobolev inequality is in force, see \cite[Section 2.4]{LedouxConcentrationLogSobolev}, or \cite{RothausConcentration}. Thus \cref{thm:NonAbelianHtype}, together with a standard argument, gives the following exponential integrability estimate similarly to Fernique's theorem for Gaussian measures. Observe that, interestingly, even if one cannot choose $\vartheta=4$ in \eqref{eqn:LogSobolevTypeUpToAdditive}, we get the upper bound $\alpha<1/(4t)$ in \cref{cor:GaussianConcentration}, which is the same as in the Euclidean space. In order to obtain the latter it is essential that we have proved that $4=\inf \mathrm{DLS}(\mathbb G)$ in \cref{thm:NonAbelianHtype}.

\begin{corollary}[Fernique's Theorem]\label{cor:GaussianConcentration}
    Let $\mathbb G$ be an $H$-type Carnot group endowed with the subRiemannian distance whose first layer has topological dimension $2n$, and whose second layer has topological dimension $m$. Let us fix $t>0$. For every $0<\alpha<\frac{1}{4t}$, and $\beta>0$, there exists a constant $C:=C(\alpha,\beta,m,n,t)$ such that the following holds.

    For every $1$-Lipschitz function $g:\mathbb G\to\mathbb R$ with $\int_{\mathbb G}g^2\de\mu_t\leqslant\beta$, we have
    \[
    \int_{\mathbb G} e^{\alpha g^2}\de\mu_t \leqslant C.
    \]
\end{corollary}

The proof of \cref{cor:GaussianConcentration} is presented in \cref{sec:Concentration}. Let us further remark that it is a classical fact that a defective log-Sobolev inequality is equivalent to  hypercontractivity of the heat semigroup, see, e.g., Gross' original paper \cite{LSIGross}, or \cite[Equations (7.5) and (7.6)]{LedouxConcentrationLogSobolev}. Hence hypercontractivity follows from our result in Theorem \ref{thm:Main1}.

Finally, in \cref{sec:defectivelogsobolevwithdistance} we exploit the same strategy used in the proof of \cref{thm:NonAbelianHtype} in order to prove a defective log-Sobolev-type inequality for the Gaussian-like measure with respect to the sub-Riemannian distance in any sub-Riemannian $H$-type group, see \cref{thm:Main1}. Similar results can be obtained in all groups satisfying a $\operatorname{MCP}$ condition, see \cref{rem:GeneralizingMainTheorem2}.
\medskip

\begin{acknowledgement} The first author acknowledges the financial support of the Courant Institute, and the AMS-Simons Travel Grant. He is also grateful to Federico Glaudo for discussions on how to compute $\min W_{1,\vartheta}$ in \eqref{eqn:Wtomega} in a computer-assisted way.
The second author is a member  of the 	Gruppo Nazionale per l'Analisi Matematica, la Probabilit\`a e le	loro Applicazioni (GNAMPA) of the Istituto Nazionale di Alta Matematica (INdAM), and  was partially funded by the INdAM-GNAMPA Project CUP\_E53C22001930001.
The third author's research was supported in part by NSF Grant DMS-2246549.
\end{acknowledgement}

\section{Preliminaries}

\subsection{H-type groups and hypoelliptic heat kernels}\label{HtypeGroups}

We recall some basic facts about $H$-type groups. For basic definitions and theory about stratified (Carnot) groups, we refer the reader to \cite{BonfiglioliLanconelliUguzzoniBook}, \cite{ledonneprimer}. An $H$-type group $\mathbb G$ is a step-2 stratified group whose Lie algebra $\mathfrak{g}$ is endowed with an inner product $\langle\cdot,\cdot\rangle$ such that
\begin{enumerate}
    \item\label{item1Htype} If $\mathfrak{z}$ is the centre of $\mathfrak{g}$ and $\mathfrak{v} = \mathfrak{z}^\perp$, then $[\mathfrak{v}, \mathfrak{v}] = \mathfrak{z}$;
\item\label{item2Htype} for every $Z \in \mathfrak{z}$, the map $J_Z : \mathfrak{v} \to \mathfrak{v}$ defined by
\[
\langle J_Z X, Y \rangle := \langle Z,[X, Y ]\rangle, \qquad \text{ for any }  X, Y \in\mathfrak{v},
\]
is an isometry whenever $\langle Z,Z\rangle =1$
\end{enumerate}

We consider the stratification $\mathfrak g=\mathfrak v \oplus \mathfrak z$. We can identify $\mathbb G$ with $\mathbb R^{2n}\times \mathbb R^m$ by means of the exponential coordinates with respect to an orthonormal basis of left-invariant vector fields. Denote the corresponding coordinates by $(x,z)$ with $x\in \mathbb R^{2n}$ and $z\in\mathbb R^m$. There is a family of anisotropic dilations $\{\delta_\lambda\}_{\lambda>0}$ defined as $\delta_\lambda(x+z):=\lambda x+\lambda^2z$, where $x\in \mathfrak v, z\in \mathfrak z$. In these coordinates we can write explicit expressions (in terms of the matrices representing $J_{\partial_{z_i}}$) of the basis of left-invariant vector fields $\mathcal{B}:=\{X_1,\ldots,X_{2n},Z_1,\ldots,Z_{m}\}$, extending the frame $\{\partial_{x_1},\ldots,\partial_{x_{2n}},\partial_{z_1},\ldots,\partial_{z_m}\}$ at the origin: see, e.g., \cite[Chapter 18]{BonfiglioliLanconelliUguzzoniBook}. 

Let $\nabla$ and $\Delta$ be the  horizontal gradient and   sub-Laplacian related to the horizontal basis $\mathcal{B}':=\{X_1,\ldots,X_{2n}\}$. Recall $\nabla f=\sum_{i=1}^{2n}(X_jf)X_j$, and $\Delta f=\sum_{i=1}^{2n}X_i^2f$.
We stress that $\Delta$ is independent of the choice of the orthonormal basis $\mathcal{B}'$, see \cite[Theorem 3.6]{GordinaLaetsch2016}. For a discussion on different notions of sub-Laplacians in sub-Riemannian manifolds, see \cite{GordinaLaetsch2016}.

In the coordinates $g=(x,z)$ introduced above $\mathscr{L}^{2n+m}$ is a (bi-invariant) Haar measure.
Let $\de\mu_t(g,\de g')$ be the transition kernel associated to the Markov semigroup $(P_t)_{t\geqslant 0}$ whose infinitesimal generator is $\Delta$.  It is known that there exists $p_t:\mathbb G\times\mathbb G\to\mathbb [0,+\infty)$ such that 
\[
\de\mu_t(g,\de g')=p_t(g,g')\de g'.
\]
Moreover, the heat kernel $p_t$ is left-invariant, i.e., $p_t(g,g')=p_t(0,g^{-1}\cdot g')$. We define, with an abuse of notation, the (hypoelliptic) heat kernel $p_t:\mathbb G\to [0,+\infty)$ as 
\[
p_t(x,z):=p_t(0,(x,z)).
\]

Recall that by \cite[Theorem 3.1(i)]{Folland1975a} one has $p_{t\varepsilon^{-2}}(x,z)=\varepsilon^{2(m+n)}p_t(\varepsilon x,\varepsilon^2z)$, and more generally such space-time scaling property holds on homogeneous Carnot groups by \cite[Remark~5.1]{CarfagniniGordina2024}.

A distinguished class of examples of $H$-type groups is given by the Heisenberg groups. Let $\mathbb H^n$ denote the (isotropic) $n$-th Heisenberg group $\mathbb H^n$, with $n\geq 1$. In this case the first layer has dimension $2n$, and the second layer has dimension $1$. We identify $\mathbb H^n$ with $\mathbb{R}^{2n+1}$ as a homogeneous Carnot group using the (global) exponential map,  with coordinates $(x_1,\ldots, x_n, y_1, \ldots, y_n, z)$. A basis of left-invariant vector fields is given by
\begin{equation}\label{eqn:HorizontalHn}
X_i=\partial_{x_i}+\frac{y_i}{2}\partial_z, \quad Y_i=\partial_{y_i}-\frac{x_i}{2}\partial_z, \qquad \text{ for any } i=1,\ldots,n,
\end{equation}
\begin{equation}
Z=\partial_z.
\end{equation}

\subsection{Sobolev inequalities on $H$-type groups}
We will use the following result which has been very recently proved in \cite[Theorem 1.2]{YangQ2024}. It is a generalization to $H$-type groups of a celebrated result by Jerison--Lee proved for $\mathbb H^n$ in \cite{JerisonLee1988}.
\begin{theorem}[\cite{JerisonLee1988, YangQ2024}]\label{thm:JerisonLee}
For every $m,n\geqslant 1$, there exists a constant $C_{n,m}$ such that the following holds. For every $H$-type group $\mathbb G$ with $\mathrm{dim}(\mathfrak{v})=2n$, and $\mathrm{dim}(\mathfrak{z})=m$, and for every $f\in W^{1,2}(\mathbb G,\mathrm{d}\mu)$  we have
\begin{equation}\label{eqn:JerisonLee}
\left(\int_{\mathbb G}|f|^{\frac{2Q}{Q-2}} \de\mu\right)^{\frac{Q-2}{Q}}\leqslant C_{n,m}\int_{\mathbb G}|\nabla f|^2\de\mu,
\end{equation}
where $Q:=2n+2m$. 
\end{theorem}

\begin{remark}\label{rem:JerisonLee}
The value of the sharp constant $C_{n,m}$ and the expression of the minimizers for arbitrary $H$-type groups are computed explicitly in \cite[Theorem 1.2]{YangQ2024}. Namely,
\[
C_{n,m}\coloneqq \frac{4^{2m/Q}}{2 n (Q-2) \pi^{(2n+m)/Q}} \left( \frac{\Gamma(2n+m)}{\Gamma(\frac{2n+m}{2})}\right)^{1/Q}.
\]
\end{remark}

As a corollary of the previous estimate, following Beckner--Pearson \cite{BecknerPearson1998}, we can deduce the following. Compare with the recent \cite[Corollary 3.6]{ChatzakouKassymovRuzhansky2023}. For discussions on the optimality of the constant $C_{n,m}$ in \cref{cor:BecknerPearson}, see \cref{QuestionLogSobolevFormLog}. 

\begin{corollary}\label{cor:BecknerPearson}
Let $\mathbb G$ be an $H$-type group. Let $C_{n,m}$ be the constant in  \eqref{eqn:JerisonLee} in \cref{thm:JerisonLee}. For every $f\in W^{1,2}(\mathbb G,\mathrm{d}\mu)$ with $\int_{\mathbb G} f^2\de\mu = 1$ we have 
\begin{equation}\label{eqn:LogSobolevForm2}
\int_{\mathbb G} f^2\log(f^2)\de\mu \leqslant (n+m)\log\left(C_{n,m}\int_{\mathbb G}|\nabla f|^2\de\mu\right).
\end{equation}
Moreover, if there exist $C\in\mathbb R$ and $D>0$ such that for every $f\in W^{1,2}(\mathbb G,\mathrm{d}\mu)$ with $\int_{\mathbb G} f^2\de\mu = 1$ we have 
\begin{equation}\label{eqn:ClogD}
\int_{\mathbb G} f^2\log(f^2)\de\mu \leqslant C\log\left(D\int_{\mathbb G}|\nabla f|^2\de\mu\right),
\end{equation}
then $C= n+m$.
\end{corollary}
\begin{proof}
Let $Q:=2m+2n$, and let $p:=2Q/(Q-2)$. Then, using Jensen's inequality
and \cref{thm:JerisonLee}, respectively, we have
\begin{equation}
\begin{split}
\frac{p-2}{2}\int_{\mathbb G} \log(|f|^2)|f|^2\de\mu &= \int_{\mathbb G} \log(|f|^{p-2})|f|^2\de\mu \leqslant \log\left(\int_{\mathbb G} |f|^p\de\mu\right)\\
&\leqslant \frac{p}{2}\log\left(C_{n,m}\int_{\mathbb G}|\nabla f|^2\de\mu\right),
\end{split}
\end{equation}
and finally observe that $p/(p-2)=n+m$ since $Q=2n+2m$.

{
The fact that $(n+m)$ is the only possible $C$ in \eqref{eqn:ClogD} comes from a scaling argument. Calling $f_\lambda(x):=\lambda^{Q/2}f(\delta_\lambda(x))$ we get that, for every $\lambda>0$, $\int_{\mathbb G}f^2\de\mu =1$ if and only if $\int_{\mathbb G} f_\lambda^2\de\mu=1$. Moreover, for every $f$ such that $\int_{\mathbb G} f^2\de\mu=1$, we have 
\[
\int_{\mathbb G} f_\lambda^2\log(f_\lambda^2)\de\mu= Q\log\lambda + \int_{\mathbb G} f^2\log(f^2)\de\mu,
\]
and
\[
\int_{\mathbb G} |\nabla f_\lambda|^2\de\mu = \lambda^2\int_{\mathbb G} |\nabla f|^2\de\mu.
\]
Thus, the validity of \eqref{eqn:ClogD} for every $f_\lambda$ with $\lambda>0$ forces $2C-Q=0$, from which the assertion follows.
}
\end{proof}

We stress here the following question, which, taking into account \cref{lem:HowToObtainLogSobolevType}, appears to be related to finding the best constant in \eqref{eqn:LogSobolevClassica}.
\begin{question}\label{QuestionLogSobolevFormLog}
Find the best (smallest) constant $\alpha_n$ for which the following holds. 
For every $f\in W^{1,2}(\mathbb H^n,\mathrm{d}\mu)$ with $\int_{\mathbb H^n} f^2\de\mu = 1$ we have 
\begin{equation}\label{eqn:Question}
\int_{\mathbb H^n} f^2\log(f^2)\de\mu \leqslant (n+1)\log\left(\alpha_n\int_{\mathbb H^n}|\nabla f|^2\de\mu\right).
\end{equation}
\end{question}

\begin{remark}\label{rem:alphanoptimal}
For the Euclidean version of \cref{QuestionLogSobolevFormLog} see \cite{Del_PinoDolbeault2003}, and the recent  \cite{BaloghDonKristaly2024}. In the Euclidean case, \eqref{eqn:Question} is equivalent to Gross' logarithmic inequality, so that optimality is achieved by replacing $n+1$ with $n/2$, the sharp constant being $\alpha_n:=2/(\pi n e)$.
    
In the Heisenberg group a version of \eqref{eqn:Question} has been obtained in \cite[Corollary 2.5]{FrankLieb2012} as a consequence of differentiating a sharp Hardy--Littlewood--Sobolev inequality on $\mathbb H^n$. Anyway, to the best of our knowledge, it is not known whether \cite[Corollary 2.5]{FrankLieb2012} gives information on the optimal constant $\alpha_n$ above. 
\end{remark}

We will need the following elementary inequality. Let $\alpha,\beta>0$ and $\gamma\in\mathbb R$, and consider the map
\begin{equation}
\delta \longmapsto \alpha\delta^2-\beta\log\delta+\gamma=:f(\delta), \delta \in (0,\infty).
\end{equation}
Then $f(\cdot)$ has a global minimum attained at $\delta_0:=\left(\beta/(2\alpha)\right)^{1/2}$. The minimal value is equal to
\begin{equation}
f(\delta_0)=\frac{\beta}{2}\log\left(e^{1+2\gamma/\beta}\cdot\frac{2\alpha}{\beta}\right).
\end{equation}
Thus for every $\alpha,\beta>0$ and $\gamma\in\mathbb R$ we get the inequality 
\begin{equation}\label{eqn:InequalityAlphaBetaGamma}
f(\delta_0)=\frac{\beta}{2}\log\left(e^{1+2\gamma/\beta}\cdot\frac{2\alpha}{\beta}\right)\leqslant \alpha+\gamma=f(1).
\end{equation}

We are now ready to prove a consequence of \cref{cor:BecknerPearson}. 

\begin{corollary}\label{cor:AnotherInequality}
Let $\mathbb G$ be an $H$-type group. Let $\tau>0$ be fixed. Let $C_{n,m}$ be the constant in \cref{thm:JerisonLee}. For every $f\in W^{1,2}(\mathbb G,\mathrm{d}\mu)$ with $\int_{\mathbb G} f^2\de\mu = 1$ we have 
\begin{equation}\label{eqn:StillAnotherInequality}
\int_{\mathbb G} f^2\log(f^2)\de\mu \leqslant \tau\int_{\mathbb G} |\nabla f|^2\de\mu + (m+n)\log\left(\frac{(m+n)C_{n,m}}{e\tau}\right).
\end{equation}
\end{corollary}
\begin{proof}
The result follows from \cref{cor:BecknerPearson} and \eqref{eqn:InequalityAlphaBetaGamma} with the choices
\begin{align*}
& \alpha:=\tau\int_{\mathbb G}|\nabla f|^2, 
\\
& \beta:=2(n+m), 
\\
& \gamma:=(m+n)\log\left(\frac{(m+n)C_{n,m}}{e\tau}\right).\qedhere
\end{align*}
\end{proof}

\subsection{BE inequality and log-Sobolev inequality in CD spaces}\label{sec:logsobolevCD}
To put our \cref{thm:NonAbelianHtype} into perspective, we register two results that are obtained by joining together several known techniques. The reader might consult \cite{BakryGentilLedouxBook} for the relevant definitions.

Let $(X,\mathcal{F},\mu)$ be a measure space, and let $(P_t)_{t\geqslant 0}$ be a Markov semigroup on it. Assume that the infinitesimal generator $\Delta$ of $(P_t)_{t\geqslant 0}$ is a diffusion operator. Let $\mathcal{A}_0$ be a \emph{standard algebra of functions} (that contains the constants when $\mu$ is a probability measure), and let $\Gamma$ be the carré du champ operator for $\Delta$ defined on $\mathcal{A}_0$ by using integration by parts, then we have that $(X,\mu,\Gamma)$ is a \emph{standard Markov triple}. See \cite[Section 3.4.2, Section 3.4.5]{BakryGentilLedouxBook} for the properties satisfied by $(X,\mu,\Gamma)$. Note that when we work in the setting of Carnot groups, we always have a Dirichlet form induced by the horizontal gradient, see \cite[p.~1902]{CarfagniniGordina2024}. Then we are automatically in the setting of a Markov triple.

We refer to \cite[Definition 3.3.14]{BakryGentilLedouxBook} for what it means for a standard Markov triple to satisfy a  curvature-dimension condition $\mathrm{CD}(\rho,n)$, where $\rho\in \mathbb R$, and $n\geqslant 1$.
\begin{proposition}\label{prop:C>2}
Let $M$ be a smooth manifold endowed with a measure $\mu$. Let $(P_t)_{t\geqslant 0}$ be a Markov semigroup, with associated diffusion operator $\Delta$, and associated carré du champ $\Gamma$. Assume that $(M,\mu,\Gamma)$ is a standard Markov triple as recalled above.

    \begin{enumerate}
        \item If there are constants $C,\varepsilon>0$ such that for every $t\in(0,\varepsilon)$ and every $f\in\mathcal{A}_0$ we have
\begin{equation}\label{eqn:1BeC}
\sqrt{\Gamma(P_tf,P_tf)}\leqslant CP_t\sqrt{\Gamma(f,f)},
\end{equation}
then $C\geqslant 1$. If $(M,\mu,\Gamma)$ does not satisfy the $\operatorname{CD}(0,\infty)$ condition, then $C>1$.
\item If there are constants $K,\varepsilon>0$ such that for every $t\in (0,\varepsilon)$, and $f \in \mathcal{A}_0$ we have
\begin{equation}\label{eqn:LSIKy}
    P_t(f^2\log f^2)-P_tf^2\log(P_tf^2) \leqslant KtP_t\Gamma( f,f),
    \end{equation}
    then $K\geqslant 4$. If $(M,\mu,\Gamma)$ does not satisfy the $\operatorname{CD}(0,\infty)$ condition, then $K>4$.
    \end{enumerate}
\end{proposition}
\begin{proof}
Since the content of this proposition is standard, we only sketch the proof for the reader's convenience. The computations in \cite[Section 5.5.1]{BakryGentilLedouxBook} give that, if \eqref{eqn:1BeC} holds, then for every $t\in (0,\varepsilon)$, and every $f\in\mathcal{A}_0$ we have
\begin{equation}\label{eqn:LSIC^2y}
P_t(f^2\log f^2)-P_tf^2\log(P_tf^2) \leqslant 4C^2tP_t\Gamma(f,f).
\end{equation}
This is~\eqref{eqn:LSIKy} with $K=4C^2$. Thus, (1) is a consequence of (2) and we may reduce to proving the latter. 
By taking $f\mapsto 1+\eta g$ in \eqref{eqn:LSIKy}, and then $\eta\to 0$, we get that for every $t\in(0,\varepsilon)$, and every $g\in\mathcal{A}_0$
\begin{equation}\label{eqn:LSIC^2ynew}
P_t(g^2)-(P_tg)^2\leqslant KtP_t\Gamma(g,g),
\end{equation}
compare with \cite[Proposition 5.1.3]{BakryGentilLedouxBook}.
Now observe that $P_th=h+t\Delta h + \frac{t^2}{2}\Delta\Delta h + o(t^2)$, and expand around $t=0$ in \eqref{eqn:LSIKy} up to the second order (compare with \cite[Theorem 4.7.2]{BakryGentilLedouxBook}). Looking at the first order one gets $K/2\geqslant 2$, from which $K\geqslant 4$. 
Finally, if $K=4$, looking at the second order, one gets that $(M,\mu,\Gamma)$ satisfies the $\mathrm{CD}(0,\infty)$ condition. Hence, if $(M,\mu,\Gamma)$ does not satisfy the $\mathrm{CD}(0,\infty)$ condition, we have $K>4$.
\end{proof}

We stress that the forthcoming \cref{thm:NonAbelianHtypePre} puts together various classical facts. The inequality in item (1) is known as the Driver--Melcher inequality \cite{DriverMelcher2005} for $p>1$, or the H.-Q.~Li inequality \cite{LiHQ2006}, see e.g. \cite{BakryBaudoinBonnefontChafai2008}. The forthcoming item (1) at that level of generality is one of the results in \cite{EldredgePhDThesis}.

\begin{theorem}\label{thm:NonAbelianHtypePre}
Let $\mathbb G$ be an H-type Carnot group and  $\de\mu$ a Haar measure on $\mathbb G$. Let $(P_t)_{t\geqslant 0}$ be the Markov semigroup whose infinitesimal generator is the sub-Laplacian $\Delta$ associated to $\langle\cdot,\cdot\rangle$. Denote by $\nabla$ the horizontal gradient in an orthonormal basis.
\begin{enumerate}
\item Let $\widetilde K:=\widetilde K(\mathbb G)$ be the minimum constant $K>0$ for which 
\begin{equation}\label{eqn:1BEinequality}
|\nabla P_t f|\leqslant K P_t|\nabla f|
\end{equation}
holds for every $t>0$, and every $f\in C^{\infty}_{\mathrm{c}}(\mathbb G)$. Then $\widetilde K\geqslant  \sqrt{\frac{3n+5}{3n+1}}$, $n$ being half of the dimension of the first layer
of $\mathbb G$.

\item Let $\mu_t$ be the heat kernel measure associated to the Markov semigroup $(P_t)_{t\geqslant 0}$. There exists a constant $C>0$ such that for every $t>0$, and every $f\in C^{\infty}_c(\mathbb G)$ with $\int f^2 \de\mu_t=1$, we have
\begin{equation}\label{eqn:LogSobolevHtype}
\int_{\mathbb G} f^2\log (f^2) \de\mu_t\leqslant Ct\int_{\mathbb G} |\nabla f|^2  \de\mu_t.
\end{equation}
    
Moreover, let $\widetilde C:=\widetilde C(\mathbb G)$ be the minimum constant $C$ for which \eqref{eqn:LogSobolevHtype} holds for every $t>0$, and every $f\in C^{\infty}_c(\mathbb G)$ with $\int f^2 \de\mu_t=1$. Then $4\widetilde K^2\geqslant \widetilde C>4$.
\end{enumerate}
\end{theorem}

\begin{proof}
 The first item directly follows from \cite[Chapter 5]{EldredgePhDThesis}, generalizing previous results in \cite{DriverMelcher2005} and \cite{LiHQ2006, LiHQ2007}. 
    
Let $\widetilde K$ be the optimal constant for which \eqref{eqn:1BEinequality} holds. By \eqref{eqn:LSIC^2y}, we get that the inequality \eqref{eqn:LogSobolevHtype} holds with $C=4\widetilde K^2$.
This proves the inequality $4\widetilde K^2\geqslant \widetilde C$. Now observe that $(\mathbb G,\de\mu,\nabla)$ does not satisfy the $\mathrm{CD}(0,\infty)$ condition by \cite[Proposition 3.6]{AmbrosioStefaniG2020}, see also \cite{RizziStefaniG2023} for recent developments. Hence, by combining the latter observation with the second part of \cref{prop:C>2}, we get that $\widetilde C> 4$, concluding the proof.
\end{proof}

\section{Proof of \cref{thm:NonAbelianHtype}}\label{proof2}
We start by sketching a simple Lemma. A result similar to \cref{lem:HowToObtainLogSobolevType} has appeared recently in \cite[Theorem 4.4]{ChatzakouKassymovRuzhansky2023}.
\begin{lemma}\label{lem:HowToObtainLogSobolevType}
Let $C, D\in\mathbb R$, let $\mathbb G$ be an $H$-type group with a Haar measure $\mu$, and $\rho> 0$ such that $\rho^{1/2},\rho^{-1/2}\in C^\infty(\mathbb G)$. Let $\xi:=\log(\rho)$ and assume that 
$\xi \in W^{1,2}_{\mathrm{loc}}(\mathbb G)$ and $\Delta\xi$ is a signed Radon measure.
The following statements are equivalent.
\begin{enumerate}
\item For every $g\in C^\infty_c(\mathbb G)$ with $\int_X g^2\de \mu=1$ we have
\begin{equation}\label{eqn:PreLogSobolevType}
\int_{\mathbb G} g^2\log (g^2) \de \mu\leqslant D+ C\int_{\mathbb G} |\nabla g|^2 \de \mu+\int_{\mathbb G} g^2\left(\frac{C}{4}|\nabla \xi|^2+\xi\right)\de \mu + \frac{C}{2}\int_{\mathbb G} g^2\de(\Delta\xi).
\end{equation}
\item For every $f\in C^\infty_c(\mathbb G)$ with $\int_{\mathbb G} f^2 \rho\de\mu=1$ we have
\begin{equation}\label{eqn:LogSobolevType}
\int_{\mathbb G} f^2\log (f^2) \rho\de\mu\leqslant D+ C\int_{\mathbb G} |\nabla f|^2  \rho\de\mu.
\end{equation}
\end{enumerate}
\end{lemma}

\begin{proof}
Let us sketch the implication (1) $\Rightarrow$ (2), the other being analogous. Suppose $f\in C^\infty_c(\mathbb G)$ with $\int_{\mathbb G} f^2\rho\de\mu=1$. If we take $g:=f\rho^{1/2}$, then $g\in C^\infty_c(\mathbb G)$ with $\int_{\mathbb G} g^2\de\mu =1$. After simple computations using integration by parts one has that
\begin{align}
& \int_{\mathbb G} g^2\log(g^2)\de\mu = \int_{\mathbb G} f^2\log(f^2)\rho\de\mu + \int_{\mathbb G} g^2\xi\de\mu,
\\
& \int_{\mathbb G} |\nabla g|^2\de\mu = \int_{\mathbb G} |\nabla f|^2\rho\de\mu-\frac{1}{4}\int_{\mathbb G} g^2|\nabla \xi|^2-\frac{1}{2}\int_{\mathbb G} g^2\de(\Delta \xi).
\end{align}
Thus the assertion follows.
\end{proof}

\subsection{Estimates on the heat kernel}\label{sec:EstimatesHeatKernel}
In the proof of \cref{thm:NonAbelianHtype} we need to bound from below the function $W_{t,C}$ defined in \eqref{eqn:Wtomega}. We thus need the refined estimates on (all  derivatives of) the heat kernel proved in \cite{BrunoCalzi2018}.

The explicit expression for the heat kernel $p_t$ in the coordinates $(x,z)$ described in \cref{HtypeGroups} above has been derived in \cite[Section~2.4]{EldredgePhDThesis}, and \cite{YangQZhuF2008}. Namely,
\begin{equation}\label{eqn:HeatExplicit}
p_t(x,z)=\frac{1}{(4\pi)^n (2 \pi)^m t^{n+m}} \int_{\mathfrak z} \mathrm e^{\frac i t \langle \lambda, z\rangle -\frac{\abs{x}^2}{4 t}\abs{\lambda} \coth\abs{\lambda}} \left( \frac{\abs{\lambda}}{\sinh \abs{\lambda}}\right)^n\,\mathrm d \lambda,
\end{equation}
for every $(x,z)\in \mathbb G$ and for every $t>0$.
In particular, there is an analytic function $h$ such that $p_t(x,z)= t^{-n-m} h( R/t, |z|/t )$, where $n=\dim \mathfrak v/2$, $m=\dim \mathfrak z$,  and $R=|x|^2/4$. We then define
\[
p_{t,k_1,k_2}(x,z)=t^{-n-m-k_1-k_2} \partial_1^{k_1} \partial_2^{k_2}h( R/t, |z|/t ).
\]
In addition, set
\[
\omega:=\frac{|z|}{R}, \qquad \delta:= \sqrt{\frac{R}{\pi |z|}}, \qquad \kappa:=2 \sqrt{\pi |z| R},
\]
and $\theta(\lambda)\coloneqq \frac{2 \lambda-\sin(2 \lambda)}{2 \sin^2(\lambda)}$, so that $\theta$ is a strictly increasing analytic diffeomorphism of $(-\pi,\pi)$ onto $\R$, and set $y_\omega\coloneqq \theta^{-1}(\omega)$ for every $\omega\in \R$. Then define
\[
d(x,z)\coloneqq \mathrm{d_{s R}}((x,z), (0,0))= \abs{x} \frac{y_\omega}{\sin(y_\omega)},
\]
with the obvious conventions for $\omega=0,\infty$, as in \cite[Definition 2.6]{BrunoCalzi2018}. Notice that here $\mathrm{d}_{\mathrm{sR}}$ is the sub-Riemannian distance in the group $\mathbb G$ associated to the restriction of $\langle\cdot,\cdot\rangle$ to the horizontal distribution.
\smallskip

Let us record the following estimates, see \cite[Theorem 4.2, Theorem 4.13, Theorem 4.14]{BrunoCalzi2018}. Related estimates had previously appeared in \cite[Chapter 4]{EldredgePhDThesis}; see also~\cite{LiHQ2010} for sharp estimates on the heat kernel. We write the estimates for $t=1$, and we stress that the general estimates for $p_{t,k_1,k_2}$ can be derived from the ones of $p_{1,k_1,k_2}$, together with the space-time scaling of the heat kernel recalled in \cref{HtypeGroups}. The following estimates hold for every fixed $C>1$.

\begin{itemize}
	\item if $\omega \leqslant C$ and $(x,z)\to \infty$, then
\[
p_{1,k_1,k_2}(x,z)=\frac{1}{R^{m/2}} \mathrm e^{-\frac{1}{4}d(x,z)^2} \Psi(\omega) \left( (-1)^{k_1+k_2} \frac{y_\omega^{n+k_1+k_2}\cos(y_\omega)^{k_1}}{\sin(y_\omega)^{n+k_1}}+ O\left(\frac 1 R\right) \right)
\]
for a suitable (explicit) analytic function $\Psi\colon [0,+\infty)\to (0,+\infty)$;
	
\item if $\delta\to 0^+$ and $\kappa\to +\infty$, then
\[
p_{1,k_1,k_2}(x,z)=\frac{(-1)^{k_2} \pi^{k_1+k_2}}{4^n (\pi\delta)^{n+k_1-(m+1)/2} \sqrt{2 \pi \kappa^m}} \mathrm e^{-\frac{1}{4} d(x,z)^2}\left(1+O\left(\delta+\frac 1 \kappa\right)\right);
\]
	
\item if $\delta\to 0^+$ and $\frac 1 C \leqslant \kappa \leqslant C$, then
\[
p_{1,k_1,k_2}(x,z)=\frac{(-1)^{k_2} \pi^{k_1+k_2}}{4^n (\pi\delta)^{n+k_1-(m+1)/2}  \kappa^{(m-1)/2}} \mathrm e^{-\frac{1}{4} d(x,z)^2} \mathrm e^{-\kappa} I_{n+k_1-1}(\kappa) \left(1+O\left(\delta\right)\right),
\]
where $I_{n+k_1-1}$ is the modified Bessel function of order $n+k_1-1$. Recall that for $\nu\in\mathbb Z$ and $\zeta\in\mathbb C$, we have 
\[
I_\nu(\zeta):=\sum_{k\in\mathbb N}\frac{\zeta^{2k+\nu}}{2^{2k+\nu}k!\Gamma(k+\nu+1)};
\]
 
\item if $\delta,\kappa\to 0^+$ (that is, if $z\to \infty$ and $\kappa \to 0^+$), then
\[
p_{1,k_1,k_2}(x,z)=\frac{(-1)^{k_2} \pi^{k_1+k_2}}{2^{2n+(m-1)/2}(n+k_1-1)!} |z|^{n+k_1-1-(m-1)/2} \mathrm e^{-\frac{1}{4} d(x,z)^2}\left(1+O\left(\kappa+\frac{1}{|z|}\right)\right).
\]
\end{itemize}

\begin{lemma}\label{lem:UnboundedOmega4}
Let $\mathbb G$ be an H-type Carnot group. Let $p_t$ be the hypoelliptic heat kernel associated to the sub-Laplacian $\Delta$,  and let $|\nabla u|$ denote the norm of the horizontal gradient $\nabla u$ of $u$.

For every $t>0$, let $\xi_t=\log p_{t}$, and for $C>0$ let
\begin{equation}\label{eqn:Wtomega}
W_{t,C}:=\frac{C}{4} |\nabla \xi_t|^2+\frac{C}{2} \Delta \xi_t+ \xi_t= -\frac{C}{4} \frac{|\nabla p_{t}|^2}{p_{t}^2}+\frac{C}{2} \frac{\Delta p_{t}}{p_{t}}+ \log p_{t}. 
\end{equation}

Then for every $t>0$ the following holds.
\[
\inf_{(x,z)\in\mathbb G}W_{t,C}(x,z)>-\infty \text{ if and only if } C>4t.
\]
\end{lemma}

We first prove the following auxiliary lemma.

\begin{lemma}\label{lem:W}
    Keep the hypotheses and the notation of Lemma~\ref{lem:UnboundedOmega4}. Then, for $z\neq 0$,
    \[
    W_{1,C}=-\frac{RC}{4} \frac{p_{1,1,0}^2+p_{1,0,1}^2}{p_{1,0,0}^2}+\frac{RC}{2}\frac{p_{1,2,0}+p_{1,0,2}+\frac{n}{R} p_{1,1,0}+\frac{m-1}{|z|}p_{1,0,1}}{p_{1,0,0}} + \log p_{1,0,0}.
    \]
\end{lemma}

\begin{proof}
    Consider the mapping $J\colon \mathfrak{v}\to \mathfrak{v}$ such that $\langle J_Z X, Y\rangle=\langle Z, [X,Y]\rangle$ for every $X,Y\in \mathfrak{v}$, so that $J_Z$ is an isometry when $\abs{Z}=1$. We define a trilinear mapping $A\colon \mathfrak v\times \mathfrak v\times \mathfrak z\to \R$ setting $A(X,Y,Z)\coloneqq \langle J_Z X, Y\rangle=\langle [X,Y],Z \rangle$, so that $A$ is antisymmetric in the first two variables, and $\norm{ A(X,\,\cdot\,,Z)}=\abs{X} \abs{Z}$ since $\abs{J_Z X}=\abs{X}\abs{Z}$. Denote by $e_1,\dots,e_{2n}$ an orthonormal basis of $\mathfrak v$ and by $u_1,\dots, u_m$ an orthonormal basis of $\mathfrak z$ (corresponding to $X_1,\dots, X_{2 n}$ and $Z_1,\dots, Z_m$, respectively).
    Then,
    \[
X_j=\partial_{e_j}+\frac{1}{2} A(x,e_j,\nabla_{\mathfrak z}),
\]
where $\nabla_{\mathfrak z}$ denotes the (Euclidean) gradient on $\mathfrak z$, for every $j=1,\dots, 2 n$. For simplicity, all $p_{1,k_1,k_2}$ in the sequel are assumed to be evaluated at some fixed $(x,z)\in \mathbb G$. In addition, $x$ and $z$ will be freely identified with $\sum_j x_j e_j$ and $\sum_j z_j u_j$, respectively.
Then,
\[
\partial_{e_j} p_{1,k_1,k_2} = p_{1,k_1+1,k_2} \partial_{e_j}R=\frac{x_j}{2} p_{1,k_1+1,k_2}.
\]
Analogously, for $z\neq 0$,
\[
\nabla_{\mathfrak z} p_{1,k_1,k_2}=p_{1,k_1,k_2+1}\nabla_{\mathfrak z}  \abs{z}=\sign(z) p_{1,k_1,k_2+1},
\]
where $\sign(z)=z/|z|$.
Furthermore, with some abuse of notation,
\[
 \nabla_{\mathfrak z}\sign(z) =\frac{1}{\abs{z}}\left(  I_{\mathfrak z}- \sign(z)\otimes \sign(z)    \right) ,
\]
where $\sign(z)\otimes \sign(z)$ denotes the mapping $z'\mapsto \langle z',\sign(z)\rangle \sign(z)$.
Therefore,
\[
X_j p_{1,0,0}= \frac{x_j}{2} p_{1,1,0}+\frac{1}{2} p_{1,0,1} A(x,e_j, \sign(z)),
\]
whence
\[
\abs{\nabla p_{1,0,0}}^2= R p_{1,1,0}^2+\frac{1}{4} p_{1,0,1}^2 \abs{J_{\sign(z)}x}^2= R (p_{1,1,0}^2+p_{1,0,1}^2  ).
\]
Furthermore,\footnote{$A(x,e_j,I_{\mathfrak z}- \sign(z)\otimes \sign(z)  ) $ should be interpreted as the element of $\mathfrak z$ corresponding to the linear map $z'\mapsto A(x,e_j,z') z'- A(x,e_j,\langle z',\sign(z)\rangle)\sign(z) $, that is, $\sum_{k=1}^m A(x, e_j, u_k) u_k- A(x, e_j, \sign(z))\sign(z)$.}
\[
\begin{split}
X_j^2 p_{1,0,0}&= \frac{1}{2} p_{1,1,0}+\frac{x_j^2}{4} p_{1,2,0}+\frac{1}{4} p_{1,1,1} A(x, x_j e_j,\sign(z))+\frac{1}{2} p_{1,0,1}A(e_j,e_j,\sign(z))\\
	&\quad+\frac{1}{2}A\bigg(x, e_j, \frac{x_j}{2} p_{1,1,1}\sign(z)+\frac{1}{2} p_{1,0,2} \sign(z) A(x,e_j, \sign(z))\\
	&\qquad\qquad+\frac{1}{2\abs{z}} p_{1,0,1} A(x,e_j,I_{\mathfrak z}- \sign(z)\otimes \sign(z)  ) \bigg) ,
\end{split}
\]
whence
\[
\begin{split}
\Delta p_{1,0,0}&= n p_{1,1,0} + R p_{1,2,0}+\frac{1}{4}p_{1,0,2} \abs{J_{\sign(z)} x}^2 \\
	&\quad+ \frac{p_{1,0,1}}{4 \abs{z}} \sum_{j=1}^{2 n} A(x, e_j, A(x,e_j,I_{\mathfrak z}- \sign(z)\otimes \sign(z)) ) 
\end{split}
\]
since $A(x,x,\sign(z))=A(e_j,e_j,\sign(z))=0$.
Now, as before 
\[
\frac{1}{4}p_{1,0,2} \abs{J_{\sign(z)} x}^2 = R p_{1,0,2},
\]
while
\[
\begin{split}
 \sum_{j=1}^{2 n} A(x, e_j, A(x, e_j, I_{\mathfrak z}- \sign(z)\otimes \sign(z) )) &= \sum_{j=1}^{2 n} \abs{A(x,e_j,\,\cdot\,)}^2-\sum_{j=1}^{2 n}\abs{A(x, e_j, \sign(z))}^2\\
	 &= \sum_{k=1}^m \abs{J_{u_k} x}^2 -\abs{J_{\sign(z)} x}^2\\
	 &=m\abs{x}^2-\abs{x}^2.
\end{split}
\]
Therefore,
\[
\Delta p_{1,0,0}= n p_{1,1,0} + R p_{1,2,0}+R p_{1,0,2}+  \frac{(m-1) R}{\abs{z}} p_{1,0,1} ,
\]
whence the result.
\end{proof}

\begin{proof}[Proof of Lemma~\ref{lem:UnboundedOmega4}.]
Notice that by the space-time scaling of $p_t$ we have 
\begin{equation}\label{eqn:Scaling}
W_{t,C}(x,z)=W_{1,C/t}(x/\sqrt t,z/t)-(n+m)\log t,
\end{equation}
for every $t>0$, for every $C\in \R$, and for every $(x,z)\in \mathbb G$, so that it will suffice to study $W_{1,C}$.
In addition, by Lemma~\ref{lem:W},
\[
W_{1,C}=-\frac{RC}{4} \frac{p_{1,1,0}^2+p_{1,0,1}^2}{p_{1,0,0}^2}+\frac{RC}{2}\frac{p_{1,2,0}+p_{1,0,2}+\frac{n}{R} p_{1,1,0}+\frac{m-1}{|z|}p_{1,0,1}}{p_{1,0,0}} + \log p_{1,0,0}.
\]
In addition, let us observe that, with the notation of~\cite{BrunoCalzi2018},
\[
\frac{d(x,z)^2}{4}=R+ \pi\abs{z}- \kappa q_\delta(\sigma_\delta)=R+ \pi\abs{z}-\kappa(1+ O(\delta^2))\sim \pi \abs{z}
\]
for $\delta\to 0$, thanks to~\cite[(3.8) and Lemma 3.9]{BrunoCalzi2018}, and since $R,\kappa=o(\abs{z})$ when $\delta\to 0$.

Then, in each of the four zones identified in \cref{sec:EstimatesHeatKernel} we have the following estimates:
\begin{itemize}
	\item  if $  \omega\leqslant C'$ and $(x,z)\to \infty$, then
    \[ 
    \frac{p_{1,k_1,k_2}(x,z)}{p_{1,0,0}(x,z)}= (-1)^{k_1+k_2} \frac{y_\omega^{k_1+k_2} \cos(y_\omega)^{k_1}}{\sin(y_\omega)^{k_1}}(1+O(1/R)),
    \]
    for every $k_1,k_2\in \N$, so that 
    \[
    \begin{split}
        W_{1,C}(x,z)&=\frac{RC}{4}\Big(-\frac{y_\omega^2 \cos(y_\omega)^2}{\sin(y_\omega)^2}-y_\omega^2+2  \frac{y_\omega^2 \cos(y_\omega)^2}{\sin(y_\omega)^2}+2y_\omega^2-2n\frac{y_\omega \cos(y_\omega)}{R\sin(y_\omega)}\\
            &\qquad-\frac{2(m-1) y_\omega}{\abs{z}} \Big)(1+O(1/R))-m\log \abs{x} -\frac{1}{4} d(x,z)^2 + O(1)    \\
            &=C\frac{R}{4} \frac{y_\omega^2}{\sin(y_\omega)^2} -m \log |x|-\frac 1 4 d(x,z)^2-\frac{C(m-1) y_\omega}{8 \omega}+O(1)\\
            &=(C-4)\frac{d(x,z)^2}{16}-m \log|x|+O(1),
    \end{split}
    \]
    since $y_\omega\sim \frac 3 2 \omega$ for $\omega\to 0^+$.
	In particular, if $C\leqslant 4$, then $W_{1,C}\to -\infty$ when $(x,z)\to \infty$ and $  \omega \leqslant C'$, so that we shall only consider the case $C>4$ in the following;
	
	\item if $C>4$, $\delta\to 0^+$, and $\kappa \to +\infty$, then
    \[
    \frac{p_{1,k_1,k_2}(x,z)}{p_{1,0,0}(x,z)}= \frac{(-\pi)^{k_2} }{\delta^{k_1}}(1+O(\delta+1/\kappa)),
    \]
    for every $k_1,k_2\in\N$,
    so that 
    \[
    \begin{split}
        W_{1,C}(x,z)&=\frac{RC}{4}\Big(-\frac{1}{\delta^2}-\pi^2+  \frac{2}{\delta^2}+2\pi^2+2\frac{n}{R\delta} -2\frac{(m-1)\pi}{\abs{z}}\Big)(1+O(\delta+1/\kappa))\\
            &\qquad-\Big(n-\frac{m+1}{2}\Big)\log \delta-\frac{m}{2}\log \kappa -\frac{1}{4}d(x,z)^2+O(1)\\
            &=\frac{C\pi\abs{z}}{4}-\frac{1}{4}d(x,z)^2 +o(\abs{z})\\
            &\sim(C-4) \frac{d(x,z)^2}{16},
    \end{split}
    \]
    since $\frac{1}{\abs{z}}= \frac{2\pi \delta}{\kappa}=o(1)$, $d(x,z)^2\sim 4\pi\abs{z}$, and $\log\delta,\log\kappa=o(\abs{z})$.
	
	\item if $\delta\to 0^+$ and $\frac{1}{C'}\leqslant \kappa \leqslant C'$, then
    \[
    \frac{p_{1,k_1,k_2}(x,z)}{p_{1,0,0}(x,z)}=\frac{(-\pi)^{k_2} I_{n+k_1-1}(\kappa)}{\delta^{k_1} I_{n-1}(\kappa) }(1+O(\delta)),
    \]
    so that
    \[
    \begin{split}
        W_{1,C}(x,z)&=\frac{RC }{4}\Big( -\frac{I_n(\kappa)^2}{\delta^2 I_{n-1}(\kappa)^2}-\pi^2 +2\frac{I_{n+1}(\kappa)}{\delta^2 I_{n-1}(\kappa)}+2\pi^2+\frac{2n I_n(\kappa)}{R\delta I_{n-1}(\kappa)}-\frac{2(m-1)\pi}{\abs{z} } \Big)(1+O(\delta))\\
            &\qquad-\Big(n-\frac{m+1}{2}\Big)\log \delta -\frac{1}{4}d(x,z)^2+O(1)\\
            &=\frac{C\pi \abs{z}}{4}\frac{2   I_{n-1}(\kappa)I_{n+1}(\kappa) -  I_n(\kappa)^2+ \frac{4 n  }{\kappa}I_{n-1}(\kappa )I_n(\kappa)}{I_{n-1}(\kappa)^2}-\frac{1}{4}d(x,z)^2+o(\abs{z})\\
            &= \frac{ 2 C I_{n-1}(\kappa)I_{n+1}(\kappa) -C I_n(\kappa)^2+ \frac{4 n C}{\kappa}I_{n-1}(\kappa )I_n(\kappa)   -4 I_{n-1}(\kappa)^2}{I_{n-1}(\kappa)^2}\frac{d(x,z)^2}{16} + o(d(x,z)^2),
    \end{split}
    \]
    since $R/\delta^2=\pi\abs{z}$, $R=o(1)$, $d(x,z)^2\sim 4\pi\abs{z}$, and $1/\delta= 2 \pi \abs{z}/\kappa$;
	
	\item if $C>4$ and $\delta,\kappa\to 0^+$, then
    \[
    \frac{p_{1,k_1,k_2}(x,z)}{p_{1,0,0}(x,z)}=\frac{(-1)^{k_2} (n-1)! \pi^{k_1+k_2} \abs{z}^{k_1}}{(n+k_1-1)!}(1+O(\kappa+1/\abs{z})),
    \]
    so that 
    \[
    \begin{split}
      W_{1,C}(x,z)&=\frac{RC }{4}\Big(-\frac{  \pi^2 \abs{z}^2}{n^2}-\pi^2 +2 \frac{\pi^2 \abs{z}^2}{n(n+1)}+2\pi^2+2\frac{\pi \abs{z}}{R}-2\frac{(m-1) \pi}{\abs{z}}    \Big) (1+O(\kappa+ 1/\abs{z}))\\
        &\qquad +\Big(n-1-\frac{m-1}{2}\Big) \log \abs{z}-\frac{1}{4} d(x,z)^2 + O(1)\\
        &=\frac{C }{2} \pi\abs{z}-\frac{1}{4} d(x,z)^2 + o(d(x,z)^2)\\  &\sim(2C-4)\frac{d(x,z)^2}{16},
    \end{split}
    \]
    since  $R\abs{z}=o(1)$, $d(x,z)^2\sim 4\pi\abs{z}$, and $\log \abs{z}=o(\abs{z})$.
\end{itemize}

So, if $C>4$, $W_{1,C}$ is bounded below in the first, second and fourth zone above. Let us deal with the third zone above. Observe first that $2  I_{n-1}(\kappa)I_{n+1}(\kappa) -I_n(\kappa)^2> 0$ for every $\kappa>0$, since $2  I_{n-1}(0)I_{n+1}(0) -I_n(0)^2= 0$ and
\[
\frac{d}{d\kappa}(2  I_{n-1}(\kappa)I_{n+1}(\kappa) -I_n(\kappa)^2)= I_{n-2}(\kappa)I_{n+1}(\kappa)+I_{n-1}(\kappa)I_{n+2}(\kappa)>0
\]
for every $\kappa>0$ (since $2I'_\nu=I_{\nu-1}+I_{\nu+1}$ for every $\nu\in\mathbb Z$). Consequently, in order to show that 
\[
\frac{ 2 C I_{n-1}(\kappa)I_{n+1}(\kappa) -C I_n(\kappa)^2+ \frac{4 n C}{\kappa}I_{n-1}(\kappa )I_n(\kappa)   -4 I_{n-1}(\kappa)^2}{I_{n-1}(\kappa)^2}>0
\]
for every $\kappa>0$ and for every $C>4$, it is necessary and sufficient to prove that
\[
J(\kappa)=2  I_{n-1}(\kappa)I_{n+1}(\kappa) -I_n(\kappa)^2+ \frac{4 n }{\kappa}I_{n-1}(\kappa )I_n(\kappa)   - I_{n-1}(\kappa)^2
\]
is $\geqslant0$ for every $\kappa>0$. This is true because of \cref{technicallemma}. Thus also in the third zone if $C>4$ then $W_{1,C}$ is bounded below. Thus the proof is concluded.
\end{proof}

\begin{proof}[Proof of \cref{thm:NonAbelianHtype}]

Without loss of generality we can assume $f\in C^\infty_c(\mathbb G)$ with $\int_{\mathbb G} f^2\de\mu=1$. Let $C'_{t,\vartheta}:=\min W_{t,\vartheta t}=\min W_{1,\vartheta} - (n+m)\log t>-\infty$, due to \cref{lem:UnboundedOmega4}, \eqref{eqn:Scaling}, and the fact that $\vartheta>4$. Let $\xi_t:=\log(p_t)$. By \cref{cor:BecknerPearson}, and \cref{lem:UnboundedOmega4}, we get 
\begin{align*}
\int_{\mathbb G} &f^2\log(f^2)\de\mu \leqslant \vartheta t\int_{\mathbb G} |\nabla f|^2\de\mu + (m+n)\log\left(\frac{(m+n)C_{n,m}}{e\vartheta t}\right) \leqslant
 \\ &\vartheta t\int_{\mathbb G} |\nabla f|^2\de\mu +\int_{\mathbb G} f^2\left(\frac{\vartheta t}{4}|\nabla \xi_t|^2+\xi_t\right)\de \mu + \frac{\vartheta t}{2}\int_{\mathbb G} f^2 \Delta\xi_t\de\mu + \\
 &+(m+n)\log\left(\frac{(m+n)C_{n,m}}{e\vartheta t}\right) - C'_{t,\vartheta}.
\end{align*}

Then, using \cref{lem:HowToObtainLogSobolevType}, one gets that the assertion in \cref{thm:NonAbelianHtype} is satisfied with 
\begin{equation}\label{eqn:ExplicitEta}
\begin{aligned}
\eta &:= (m+n)\log\left(\frac{(m+n)C_{n,m}}{e\vartheta t}\right) - C'_{t,\vartheta} \\
&=(m+n)\log\left(\frac{(m+n)C_{n,m}}{e\vartheta t}\right) + (n+m)\log t - \min_{g\in\mathbb G}W_{1,\vartheta}(g) \\
&=(m+n)\log\left(\frac{(m+n)C_{n,m}}{e\vartheta}\right) - \min_{g\in\mathbb G}W_{1,\vartheta}(g),
\end{aligned}
\end{equation}
    where $W_{1,\vartheta}$ is defined in \eqref{eqn:Wtomega}. Notice that in $H$-type groups there is an explicit formula for the heat kernel in the coordinates we are considering, see \eqref{eqn:HeatExplicit}. Moreover this expression only depends on $n,m$. Thus, if wanted, one might get an approximate value of $\eta$, only depending on $n,m,\vartheta$, in a computer-assisted way by estimating $\min_{g\in\mathbb G}W_{1,\vartheta}(g)$.

    Let us now show the last part of the statement. Suppose by contradiction $4\in \mathrm{DLS}(\mathbb G)$. Then, there exists a $\eta>0$ such that $\mathrm{DLS}(4,\eta)$ holds. Using \cref{lem:HowToObtainLogSobolevType}, we have that for every $f\in C^{\infty}_c(\mathbb G)$ with $\int_{\mathbb G}f^2\de\mu=1$ the following inequality holds
\begin{equation}\label{eqn:PreLogSobolevTypeNEW}
\int_{\mathbb G} f^2\log (f^2) \de \mu - 4\int_{\mathbb G} |\nabla f|^2\de \mu - \eta \leqslant\int_{\mathbb G} f^2\left(|\nabla \xi_1|^2+\xi_1+2\Delta \xi_1\right)\de \mu,
\end{equation}
where $\xi=\log(p_1)$.
Let $B^{\mathbb G}_r(p)$ denote the ball of radius $r>0$ and center $p\in \mathbb G$ with respect to the sub-Riemannian distance. Let us now fix $f$ to be a positive smooth function compactly supported in $B^{\mathbb G}_1(0)$, such that $\int_{\mathbb G}f^2\de\mu=1$, and $f(x)>c>0$ for every $x\in B^{\mathbb G}_{1/2}(0)$ and for some constant $c>0$. For every $v\in\mathbb G$ define $f_v(x):=f(v\cdot x)$, where $\cdot$ is the product in the group $\mathbb G$. Since $\mu$ is left-invariant we have that $\int_{\mathbb G} f_v^2\de \mu=1$, and there is $\Lambda\in\mathbb R$ such that
\[
\int_{\mathbb G} f_v^2\log (f_v^2) \de \mu - 4\int_{\mathbb G} |\nabla f_v|^2\de \mu - \eta = \Lambda,
\]
for every $v\in \mathbb G$. Moreover, recall that $W_{1,4}=|\nabla\xi_1|^2+\xi_1+2\Delta\xi_1$. By inspecting the asymptotics of $W_{1,4}$ in the first zone in the proof of \cref{lem:UnboundedOmega4} we get that there exists a sequence $v_i\to+\infty$ in $\mathbb G$ such that 
\begin{equation}\label{eqn:To-Inf}
\sup_{x\in B_{1}^{\mathbb G}(v_i)} \left[|\nabla\xi_1|^2(x)+\xi_1(x)+2\Delta\xi_1(x)\right] \to -\infty.
\end{equation}
Testing \eqref{eqn:PreLogSobolevTypeNEW} with $f:=f_{v_i^{-1}}$, and using that $\de\mu$ is left-invariant, we get
\[
\int_{\mathbb G} f^2\left[|\nabla\xi_1|^2(v_i\cdot x)+\xi_1(v_i\cdot x)+2\Delta\xi_1(v_i\cdot x)\right] \geqslant \Lambda,
\]
for every $i$. The latter is a contradiction with \eqref{eqn:To-Inf} and the fact that $f>c>0$ on $B_{1/2}^{\mathbb G}(0)$.
\end{proof}
\subsection{Exponential integrability of Lipschitz functions}\label{sec:Concentration}
We are now ready to prove \cref{cor:GaussianConcentration}.

\begin{proof}[Proof of \cref{cor:GaussianConcentration}]   
    Let us sketch the classical proof. Fix $\vartheta>4$ such that $1/(\vartheta t)>\alpha$. This $\vartheta$ exists and depends on $\alpha$, and $t$, since $0<\alpha<1/(4t)$. Hence, by \cref{thm:NonAbelianHtype} there exists $\eta:=\eta(\vartheta,m,n)$ such that \eqref{eqn:FullDefectiveLogSob} holds. Fix $g$ as in the statement. Notice that $|\nabla_{\mathbb G}g|\leqslant 1$ almost everywhere. By considering $\phi_n\circ g$ where $\phi_n:\mathbb R\to\mathbb R$ with $\phi_n\in C^{\infty}_c(\mathbb R)$, $|\phi_n'|\leqslant 1$, $\phi_n(x)=x $ on $(-n,n)$, $|\phi(x)|\leqslant |x|$, and eventually sending $n\to \infty$, we can assume $g$ is in addition with compact support. For every $\lambda>0$, define $K(\lambda):=\lambda^{-1}\log\left(\int_{\mathbb G}e^{\lambda g}\de\mu_t\right)$. By applying \eqref{eqn:FullDefectiveLogSob} to $f:=e^{\frac{\lambda g}{2}}\in W^{1,2}(\mathbb G,\mathrm{d}\mu_t)$ for every $\lambda>0$ we get the following ODE
    \[
    K'(\lambda)\leqslant \frac{\eta}{\lambda^2}+\frac{\vartheta t}{4}, \text{ for all } \lambda >0.
    \]
    Integrating the previous ODE from $1$ to $\lambda$ and arguing as in \cite[Proposition 2.8]{LedouxConcentrationLogSobolev} will give that there exists a constant $B:=B(\vartheta,m,n,\beta)$ such that
    \[
    \int_{\mathbb G}e^{\lambda g}\de\mu_t \leqslant e^{B\lambda+\frac{\vartheta t\lambda^2}{4}} \text{ for all } \lambda\geqslant 1,  \text{ for every $g$ such that $\int_{\mathbb G}g^2\de\mu_t\leqslant \beta$}.
    \]
    The latter, and Markov's inequality, imply that if $r\geqslant B+\vartheta t/2$ then 
    \[
    \int_{\{|g|\geqslant r\}} \de\mu_t \leqslant e^{-\frac{(r-B)^2}{\vartheta t}}.
    \]
    Finally, using Fubini (see, e.g., \cite[Proposition 1.2]{LedouxConcentrationLogSobolev}) the previous inequality implies that 
    \[
    \int_{\mathbb G}e^{\gamma g^2}\de\mu_t<C(\gamma,\vartheta,m,n,\beta,t),
    \]
    for every $\gamma < 1/(\vartheta t)$. Thus taking $\gamma=\alpha$ concludes the proof.
\end{proof}
 
\section{A defective log-Sobolev inequality  for the Gaussian-like measure with respect to the sub-Riemannian distance}\label{sec:defectivelogsobolevwithdistance}

In this section we exploit the techniques we used to prove Theorem \ref{thm:NonAbelianHtype} to show another defective log-Sobolev inequality for the Gaussian-like measure with respect to the sub-Riemannian distance in the H-type group $\mathbb G$.

We first discuss some notation. Recall that $2n$ and $m$ denote the dimensions of the first, and second layer, respectively, of the Lie algebra of $\mathbb G$. Let $\mathrm{d}\mu=\mathscr{L}^{2n+m}$ be the Lebesgue measure in the exponential coordinates associated to an orthonormal basis of the Lie algebra of $\mathbb G$, and let $\dist_{\mathrm{sR}}$ be the sub-Riemannian distance associated to the horizontal basis $\mathcal{B}\subset\mathcal{B}'$. One can show
\[
c:=\int_{\mathbb G} e^{-\dist_{\mathrm{sR}}^2(x,0)/2}\de\mu(x) < +\infty,
\]
so that $e^{-\dist_{\mathrm{sR}}^2(x,0)/2}\de\mu(x)$ is a finite measure. Notice that $e^{-\dist_{\mathrm{sR}}^2(x,0)/2}\de\mu(x)$ might not be a probability measure. By normalizing and by appropriately changing all the constants in the proof of Theorem \ref{thm:Main1} one gets a similar result with the probability measure $c^{-1}e^{-\dist_{\mathrm{sR}}^2(x,0)/2}\de\mu(x)$. While the constants we obtain are likely not to be optimal, we present Theorem \ref{thm:Main1} as follows.

\begin{theorem}\label{thm:Main1}
For every $f\in W^{1,2}(\mathbb G,e^{-\dist_{\mathrm{sR}}^2(x,0)/2}\de\mu(x))$ with 
\[
\int_{\mathbb G} f^2(x)\,(e^{-\dist_{\mathrm{sR}}^2(x,0)/2}\de\mu(x))=1,
\]
we have
\begin{equation}\label{eqn:LogSobolevTypeDistanceSquared}
\int_{\mathbb G} f^2(x)\log (f^2(x)) \, e^{-\dist_{\mathrm{sR}}^2(x,0)/2}\de\mu(x)\leqslant K_{n,m}+ 2\int_{\mathbb G} |\nabla f|^2(x) \, e^{-\dist_{\mathrm{sR}}^2(x,0)/2}\de\mu(x),
\end{equation}
where
\begin{equation}
K_{n,m}\coloneqq \left(\log\left(\left(\frac{(n+m)C_{n,m}}{2e}\right)^{n+m}\right)+2n+3m\right).
\end{equation}
Finally, if $\vartheta< 2$, then there is no $\eta\geqslant 0$ such that 
\begin{equation}\label{eqn:LogSobolevTypeDistanceSquaredConMenoDi2}
\int_{\mathbb G} f^2(x)\log (f^2(x)) \, e^{-\dist_{\mathrm{sR}}^2(x,0)/2}\de\mu(x)\leqslant \eta+ \vartheta\int_{\mathbb G} |\nabla f|^2(x) \, e^{-\dist_{\mathrm{sR}}^2(x,0)/2}\de\mu(x),
\end{equation}
for every $f\in W^{1,2}(\mathbb G,e^{-\dist_{\mathrm{sR}}^2(x,0)/2}\de\mu(x))$ with $\int_{\mathbb G} f^2(x)\,(e^{-\dist_{\mathrm{sR}}^2(x,0)/2}\de\mu(x))=1$.
\end{theorem}

The constant $K_{n,m}$ in \cref{thm:Main1} is unlikely to be optimal. Compare with the discussions in \cref{cor:BecknerPearson}, \cref{cor:AnotherInequality}, \cref{QuestionLogSobolevFormLog}, and \cref{rem:alphanoptimal}. 
Using the same technique as in the proof of \cref{thm:Main1} one can prove similar inequalities for other Carnot groups satisfying the measure contraction property, see \cref{rem:GeneralizingMainTheorem2}.

\subsection{Proof of \cref{thm:Main1}}
We start with the following lemma.
\begin{lemma}\label{lem:CBoundedBelow}
Let $\rho(x):= e^{-\mathrm{d}_{\mathrm{sR}}^2(x,0)/2}$, and $\xi(x):=\log(\rho(x))=-\mathrm{d}_{\mathrm{sR}}^2(x,0)/2$. The infimum of $C>0$ such that $\frac{C}{4}|\nabla \xi|^2+\frac{C}{2}\Delta\xi+\xi$ is bounded below (distributionally) on $\mathbb G$ is $2$.
\end{lemma}

\begin{proof}
With a slight abuse of notation, we denote $\dist_{\mathrm{sR}}(x):=\dist_{\mathrm{sR}}(x,0)$. We have 
\begin{align}
 \left|\nabla \left(-\mathrm{d}_{\mathrm{sR}}^2/2\right)\right|^2 &= \mathrm{d}_{\mathrm{sR}}^2, 
\notag
\\
 \Delta \left(-\mathrm{d}_{\mathrm{sR}}^2/2\right)&\geqslant -(2n+3m),
\end{align}
where the last inequality comes from the fact that $(\mathbb G,\dist_{\mathrm{sR}},\mu)$ satisfies  $\mathrm{MCP}(0,2n+3m)$ by \cite[Theorem 3]{BarilariRizzi2018}, and thus using \cite[Corollary 4.17]{CavallettiMondino2020a}.
Thus the assertion follows. Moreover, 
\begin{equation}\label{eqn:DistributionalBoundLaplacian}
\frac{1}{2}|\nabla \xi|^2+\Delta\xi+\xi \geqslant -(2n+3m),
\end{equation}
holds distributionally.
\end{proof}

\begin{proof}[Proof of \cref{thm:Main1}]
By density, it suffices to prove the result for $g\in C^\infty_{\mathrm{c}}(\mathbb G)$. Let $\xi:=\log(\rho)=-\mathrm{d}_{\mathrm{sR}}^2/2$. By \cref{cor:AnotherInequality} with $\tau=2$, and \eqref{eqn:DistributionalBoundLaplacian}, respectively, we have the following. For every $g\in C^\infty_{\mathrm{c}}(\mathbb G)$ with $\int_{\mathbb G} g^2\de\mu=1$ we have
\begin{align}
& \int_{\mathbb G} g^2\log(g^2)\de\mu \leqslant 2\int_{\mathbb G} |\nabla g|^2\de\mu + \log\left(\left(\frac{(n+m)C_{n,m}}{2e}\right)^{n+m}\right)
\notag
\\
&\leqslant \left(\log\left(\left(\frac{(n+m)C_{n,m}}{2e}\right)^{n+m}\right)+2n+3m\right)
\\
&+2\int_{\mathbb G} |\nabla g|^2\de\mu+\int_{\mathbb H^n}g^2\left(\frac{1}{2}|\nabla \xi|^2+\xi\right)\de\mu+\int_{\mathbb G}g^2\de(\Delta\xi).
\notag
\end{align}
Thus we can use (a weak variant of) \cref{lem:HowToObtainLogSobolevType} to complete the proof. 

The fact that the constant $2$ is optimal can be explained as follows. Suppose \eqref{eqn:LogSobolevTypeDistanceSquaredConMenoDi2} holds with a constant $\vartheta<2$, and with some $\eta\in\mathbb R$. Then the corresponding inequality \eqref{eqn:PreLogSobolevType} holds as well with $C=\vartheta<2$: indeed, this follows directly  from \cref{lem:HowToObtainLogSobolevType}. As a consequence of this, fixing an arbitrary $g\in C^\infty_{\mathrm{c}}(\mathbb G)$ with $\int_{\mathbb G} g^2\de\mu=1$, and denoting $g_v:=g(v\cdot x)$ for every $v\in\mathbb G$, we have that
\[
\int_{\mathbb G} g_v^2\left(\frac{C}{4}|\nabla\xi|+\frac{C}{2}\Delta\xi+\xi\right)\de\mu  \text{ is uniformly bounded below for $v$ big enough.}
\]
Thus, since from \cref{lem:CBoundedBelow} we get that $\left(\frac{C}{4}|\nabla\xi|+\frac{C}{2}\Delta\xi+\xi\right)$ is unbounded below when $C<2$, the previous argument gives a contradiction taking $v\to+\infty$, and $g$ appropriately.
\end{proof}

\begin{remark}\label{rem:GeneralizingMainTheorem2}
The argument used to prove \cref{thm:Main1} would give an inequality like \cref{thm:Main1}, with computable constants, for every Carnot group satisfying $\mathrm{MCP}(0, K)$ for some explicit $K>0$, and for which there is an explicit estimate of the $(2,2^{\ast})$-Sobolev constant. 
Note that \emph{an} optimal constant for the $(2,2^*)$-Sobolev constant, which is moreover attained on at least one function in $W^{1,2}$, always exists in every Carnot group, see \cite{Vassilev2006a}. To the best of our knowledge estimates of such an optimal constant in Carnot groups other than $H$-type groups are not available.

Let us further recall that an arbitrary step-$2$ Carnot group satisfies $\mathrm{MCP}(0,N)$ for some $N>1$, see \cite{BadreddineRifford}. For related results see \cite{Milman2021a}.
\end{remark}

\appendix

\section{A technical Lemma}

\begin{lemma}\label{technicallemma}
    Let $n\geqslant 1$, and let $I_{n}$ be the modified Bessel function of order $n$. 
    Then
    \[
    J(\kappa):=2  I_{n-1}(\kappa)I_{n+1}(\kappa) -I_n(\kappa)^2+ \frac{4 n }{\kappa}I_{n-1}(\kappa )I_n(\kappa)   - I_{n-1}(\kappa)^2\geqslant 0.
    \]
\end{lemma}
\begin{proof}
    Let us recall that, for every $n$, we do have \cite[7.11(23)]{HigherTranscBook1953}
    \[
    \frac{2n}{\kappa}I_n(\kappa)=I_{n-1}(\kappa)-I_{n+1}(\kappa).
    \]
    Thus
    \begin{equation}
        \begin{split}
            J(\kappa)&=2I_{n-1}(\kappa)I_{n+1}(\kappa)-I_n(\kappa)^2+2I_{n-1}(\kappa)(I_{n-1}(\kappa)-I_{n+1}(k))-I_{n-1}(\kappa)^2\\
            &=I_{n-1}(\kappa)^2-I_n(\kappa)^2\geqslant 0,
        \end{split}
    \end{equation}
    where the last inequality comes from the fact that $n\mapsto I_n$ is a pointwise decreasing family of functions.
\end{proof}

\section{More general $2$-step groups}

In this section we shall briefly comment on what our methods may achieve when $\mathbb G$ is allowed to be a more general $2$-step stratified group. 
We begin with the simplest case, that is, a Heisenberg group $\mathbb H^n$ endowed with an arbitrary inner product. In particular, in this case $\mathbb G$ need no longer be an $H$-type group, so that~\cite{BrunoCalzi2018} cannot be invoked to determine precise asymptotic estimates for the heat kernel and its derivatives. 
As a matter of fact, by means of a suitable change of coordinates, one may always reduce to the case of the non-isotropic Heisenberg groups considered, e.g., in~\cite{LiZhang, GordinaLuo2022}. In this case, the heat kernel $p_t$ may be described as follows (cf.~\cite{LiZhang}):
\[
p_t(x,z)=\frac{2}{(4 \pi t)^{n+1}} \int_\R e^{\frac 1 t \big( i z s -\sum_{j=1}^k \frac{\alpha_j \abs{P_j(x) }^2}{4} s \coth(\alpha_j s)\big)} \prod_{j=1}^k \left( \frac{\alpha_j s}{\sinh(\alpha_j s)}\right)^{n_j}\,\de s,
\]
for every $t>0$ and for every $(x,z)\in \mathbb G$, where  $1\leqslant k \leqslant n$, $n_1,\dots,n_k\geqslant 1$, $\sum_{j=1}^k n_j=n$, $0<\alpha_1<\cdots<\alpha_k$, and the $P_j$ are commuting self-adjoint projectors of the first layer $\mathfrak z^\perp$ such that $P_j P_{j'}=0$ for $j\neq j'$, $\sum_{j=1}^k P_j=I$, and  $P_j(\mathfrak z^\perp)\oplus \mathfrak z$ is isomorphic to a Heisenberg group $\mathbb H^{n_j}$.

In this situation, the method of the stationary phase may be applied (up to some extent), following essentially the procedure of~\cite{BrunoCalzi2018} with minor modifications.
In fact, one may prove that, for every $C>0$, 
\begin{equation}\label{eqn:Expansionp1anisotropic}
p_1(x,z)= \frac{2 \sqrt \pi e^{ -d(x,z)^2 /4}}{(4 \pi)^{n+1}  \Psi(x,z)}  \left[\prod_{j=1}^k\left( \frac{\alpha_j y_{(x,z)}}{\sin(y_{(x,z)})}   \right)^{n_j} + O\Big(\frac 1 R\Big)\right] 
\end{equation}
for $(x,z)\to \infty$ with $\frac{\abs{z}}{R}\leqslant C$ and $\frac{\abs{P_k(x)}^2}{R}\geqslant 1/C$, where $R=\sum_{j=1}^k \alpha_j \abs{P_j(x)}^2$,
\[
\frac 1 4 d(x,z)^2=  \abs{z} y_{(x,z)} + \sum_j \frac{\alpha_j \abs{P_j(x)}^2}{4 } y_{(x,z)} \cot(\alpha_j y_{(x,z)}) ,
\]
$y_{(x,z)}$ is the unique element of $[0,\pi/\alpha_k)$ such that
\[
\sum_j  \alpha_j \abs{P_j(x)}^2  \frac{2 \alpha_j y_{(x,z)}-\sin(2\alpha_j y_{(x,z)}) }{2 \sin(\alpha_j y_{(x,z)})^2}=4 \abs{z},
\]
and
\[
 \Psi(x,z)=\sqrt{   \sum_{j=1}^k \alpha_j^2 \abs{P_j(x)}^2 \frac{\sin(\alpha_j y_{(x,z)})-\alpha_j y_{(x,z)}\cos(\alpha_j y_{(x,z)})}{ 4\sin(\alpha_j y_{(x,z)})^3}  } .
\]
Similar estimates may be found also for the derivatives of $p_1$. Tedious computations then show that $W_{1,C'}(x,z)=(C'-4) \frac{d(x,z)^2}{16}-\log\abs{x} + O(1)$ for $(x,z)\to \infty$, $\abs{z}\leqslant C R$ and $\abs{P_k(x)}\geqslant C\sqrt R$.
When $k=1$ and $\alpha_1=1$, $\mathbb G$ is an $H$-type group and we re-obtain the previous estimates (when $\alpha_1\neq 1$ some constants differ, but the estimates are essentially the same), whereas when $k>1$ a new difficulty arises, since the above estimates tell us nothing when $\frac{\abs{P_k(x)}^2}{R}\to 0$. As a matter of fact, when $P_k(x)=0$, then one may find estimates similar to the preceding ones, replacing $k$ with $k-1$, and so on. Nonetheless, this is not enough to find ``uniform'' asymptotic estimates for $(x,z)\to \infty$ and $\frac{\abs{z}}{R}\leqslant C$. Since the picture would be incomplete in any case, we shall not examine the case $\frac{\abs{z}}{R}\to \infty$ (which is the most difficult one when $\mathbb G$ is an $H$-type group).

\begin{proof}[Proof of \eqref{eqn:Expansionp1anisotropic}]
    The idea is the same as in~\cite[Theorem 3.3]{BrunoCalzi2018}. Roughly speaking, $p_1$ is essentially written in the form
    \[
    \int_\R e^{i R f(\lambda)} g(\lambda) \,\de \lambda 
    \]
    (with $f$ and $g$ depending also on $x/\sqrt R$ and $\abs{z}/R$). The method of the stationary phase allows us to find rather precise asymptotic estimates of such an expression for $R\to +\infty$, provided that $f'(\lambda_0)=0$ and $f''(\lambda_0)\neq 0$ for some $\lambda_0\in \R$ (and some further technical assumptions are satisfied). 
    Unfortunately, there are no such $\lambda_0\in \R$ unless $z=0$; nonetheless, one may leverage the fact that $f$ and $g$ extend to holomorphic functions on a horizontal strip to find some $\lambda_0=i y_{(x,z)}$ therein with the required properties. 
    Unlike the case of H-type groups, though, it turns out that the `natural' $y_{(x,z)}$ is only defined and analytic when $P_k(x)\neq 0$; when $P_k(x)=0$, the `natural' point to choose would be larger than the limit of the $y_{(x,z)}$ for $P_k(x)\neq 0$. This issue prevents the estimates one may find with this method from extending `uniformly' to the region where $P_k(x)\to 0$. Let us then pass to (a sketch of) the proof.
    
    Let us first introduce some notation. Define
    \[
    \theta\colon z \mapsto \begin{cases}
    \frac{2 z- \sin(2 z)}{2 \sin(z)^2} &\text{if $z\not \in \pi \mathbb Z$}\\
    0 &\text{if $z=0$},
    \end{cases}
    \]
    and observe that $\theta$ is holomorphic and that 
    \[
    \theta'(z)=2\frac{\sin(z)-z \cos(z)}{\sin(z)^3}.
    \]
    In particular, $\theta$ induces a strictly increasing odd diffeomorphism of $(-\pi,\pi)$ onto $\R$. Next, define 
    \[
    \phi_{(x,z)}\colon s \mapsto \begin{cases}
        \frac{\abs{z}}{R} s +i\sum_{j=1}^k \frac{\alpha_j \abs{P_j(x) }^2}{4 R} s \coth(\alpha_j s) &\text{if $s \not \in \bigcup_j \frac{i\pi}{\alpha_j} \mathbb Z$}\\
        i &\text{if $s=0$},
    \end{cases}
    \]
    so that
    \[
    p_1(x,z)= \frac{2}{(4 \pi)^{n+1}} \int_\R e^{ i R \phi_{(x,z)}(s)} a(s)\,\de s,
    \]
    where 
    \[
    a\colon s \mapsto \begin{cases}
        \prod_j \left( \frac{\alpha_j s}{\sinh(\alpha_j s)}\right)^{n_j} &\text{if $s\not \in \bigcup_j \frac{i\pi}{\alpha_j} \mathbb Z$}\\
        1 & \text{if $s=0$}.
    \end{cases}
    \]
    Define 
    \[
    \theta_{(x,z)}\colon s \mapsto\sum_{j=1}^k \frac{\alpha_j \abs{P_j(x) }^2}{4 R} \theta( \alpha_j s),
    \]
    so that
    \[
    \phi_{(x,z)}'(s)= \frac{\abs{z}}{R} + \theta_{(x,z)}(i s).
    \]
    In addition, if $P_k(x)\neq 0$, then $\theta_{(x,z)}$ induces a strictly increasing odd analytic diffeomorphism of $(-\pi/\alpha_k,\pi/\alpha_k)$ onto $\R$, so that $y_{(x,z)}=\theta_{(x,z)}^{-1}\big(  \frac{\abs{z}}{R}\big)$ is well defined (and analytic as a function of $(x,z)$ where $P_k(x)\neq 0$). In addition, if we set
    \[
    \psi_{(x,z)}\colon s \mapsto \phi_{(x,z)}(s+ i y_{(x,z)})- \phi_{(x,z)}(i y_{(x,z)}),
    \]
    then $\psi_{(x,z)}(0)=\psi_{(x,z)}'(0)=0$ and
    \[
    \psi''_{(x,z)}(0)= i\theta'_{(x,z)}(y_{(x,z)})= i \sum_{j=1}^k \frac{\alpha_j^2 \abs{P_j(x) }^2}{2 R} \frac{\sin(\alpha_j y_{(x,z)})-\alpha_j y_{(x,z)} \cos(\alpha_j y_{(x,z)})}{\sin(\alpha_j y_{(x,z)})^3}.
    \]
    Since $i \psi''_{(x,z)}(0) < 0$, by means of Taylor's formula one may show that there are $C_1>0$ and $\eta>0$ such that 
    \[
    \abs{\psi'_{(x,z)}(s)}\geqslant C_1 \abs{s}
    \]
    for every $(x,z)$ such that $\frac{  \abs{P_k(x) }^2}{R}\geqslant 1/C$ and $\frac{\abs{z}}{R}\leqslant C$, and for every $s\in [-\eta,\eta]$.

    Now, observe that~\cite[Lemma 3.4]{BrunoCalzi2018} shows that 
    \[
    \mathrm{Im} \phi_{(x,z)}(s) > \frac{\abs{z}}{R} \mathrm{Im} s
    \]
    provided that $\abs{\mathrm{Re} s}>\abs{\mathrm{Im} s}$, and that
    \[
    \abs{a(s)}= \abs{s}^n \prod_j \frac{\alpha_j^{n_j}}{(\sinh(\alpha_j \mathrm{Re} s)^2+ \sin(\alpha_j \mathrm{Im} s)^2)^{n_j/2}}
    \]
    when $s\not \in \bigcup_j \frac{i\pi}{\alpha_j} \mathbb Z$. Therefore, we may change the contour of integration and show that
    \[
    p_1(x,z)= \frac{2 e^{i R \phi_{(x,z)}(i y_{(x,z)})}}{(4 \pi)^{n+1}} \int_\R e^{i R \psi_{(x,z)}(s)} a(s+ i y_{(x,z)}) \,\de s.
    \]
    Observe that $i R \phi_{(x,z)}(i y_{(x,z)})=-d(x,z)^2/4$ by definition.
    Next, observe that $t \cot t\leqslant 1$ for every $t \in (-\pi,\pi)$, so that
    \[
    \begin{split}
    \mathrm{Im} \psi_{(x,z)}(s)&=\sum_j \frac{\alpha_j \abs{P_j(x) }^2}{4 R} \sinh(\alpha_j s)^2\frac{ \alpha_j s \coth(\alpha_j s)  - \alpha_j y_{(x,z)} \cot(\alpha_j y_{(x,z)})  }{\sinh(\alpha_j s)^2+\sin(\alpha_j y_{(x,z)})^2}\\
        &\geqslant \sum_j \frac{\alpha_j \abs{P_j(x) }^2}{4 R} \frac{\alpha_j s \coth(\alpha_j s) -1}{1+\sinh(\alpha_j s)^{-2}}>0
    \end{split}
    \]
    and that $\frac{\alpha_j s \coth(\alpha_j s) -1}{1+\sinh(\alpha_j s)^{-2}}\sim \abs{s} $ for $s\to \infty$.
    Therefore, the method of the stationary phase (cf., e.g.,~\cite[Theorem 2.7]{BrunoCalzi2018}) allows to conclude that
    \[
    p_1(x,z) =  \frac{2 e^{-d(x,z)^2/4}}{(4 \pi)^{n+1}} \left[\sqrt{ \frac{2 \pi i}{R \phi''_{(x,z)}(0)} } a(i y_{(x,z)}) + O(1/R^{3/2})\right]
    \]
    whence our claim.    
\end{proof}

In the general case, the situation is even more intricate. Indeed, one may determine an expression of $p_t$ somewhat similar to the previous one (cf., e.g.,~\cite[Theorem 4.23]{YangZPhDThesis2022}). In this case, however, the $\alpha_j$ and the $P_j$ \emph{depend} on $s$, and may be defined as the `eigenvalues' and the `eigenprojectors' (in a sense) relative to the diagonalization of the inner product with respect to the skew-symmetric form $(x,x')\mapsto \langle [x,x'] , s \rangle$ (when this form is not symplectic, one has to add a Gaussian factor $e^{-c \abs{P_0(x)}^2}$, where $P_0$ (depends on $s$ and) is the self-adjoint projector onto the radical of the preceding skew-symmetric form). Notice that, even though the $\alpha_j$ and the $P_j$ are generically well defined and real analytic, they may behave badly on an algebraic set; in particular, the $\alpha_j$ may always be chosen so as to be (H\"older) continuous, but need not be even Lipschitz. Since the method of the stationary phase may be applied to this kind of integral expression only after a translation to a line of the form $\R+ i y_{(x,z)}$, this lack of analyticity introduces a severe obstacle to a fruitful extension of the methods used in~\cite{BrunoCalzi2018}.

\printbibliography

\end{document}